\documentclass{amsart}
\usepackage{graphicx,amssymb,epsfig}
\usepackage[all]{xy}

\newtheorem{theorem}{Theorem}[section]
\newtheorem{proposition}[theorem]{Proposition}
\newtheorem{lemma}[theorem]{Lemma}
\newtheorem{corollary}[theorem]{Corollary}
\newtheorem{remark}[theorem]{Remark}
\newtheorem{example}[theorem]{Example}

\newtheorem{definition}[theorem]{Definition}

\newcommand{\bth}{\begin{theorem}}
\newcommand{\bpr}{\begin{proposition}}
\newcommand{\epr}{\end{proposition}}
\newcommand{\bco}{\begin{corollary}}
\newcommand{\eco}{\end{corollary}}
\newcommand{\bde}{\begin{definition}\rm}
\newcommand{\ede}{\end{definition}}
\newcommand{\ble}{\begin{lemma}}
\newcommand{\ele}{\end{lemma}}
\newcommand{\bre}{\begin{remark}\rm}
\newcommand{\ere}{\end{remark}}
\newcommand{\bex}{\begin{example}\rm}
\newcommand{\eex}{\end{example}}

\def\la#1{\hbox to #1pc{\leftarrowfill}}
\def\ra#1{\hbox to #1pc{\rightarrowfill}}
\def\ma#1{\hbox to #1pc{\mapstofill}}

\def\fract#1#2{\raise4pt\hbox{$ #1 \atop #2 $}}

\def\lrar{{\ra 2}}

\def\tensor{\otimes}
\def\sp#1{\hbox{SP}^{#1}}

\def\wsp#1{\widehat{\hbox{SP}}^{#1}}
\def\bsp#1{\overline{\hbox{SP}}^{#1}}

\def\sub#1{\hbox{Sub}_{#1}}
\def\bsub#1{\overline{\hbox{Sub}}_{#1}}
\def\sn{{\mathfrak S_n}}

\def\bbz{{\mathbb Z}}
\def\bbf{{\mathbb F}}
\def\bbp{{\mathbb P}}
\def\bbr{{\mathbb R}}

\def\bbc{{\mathbb C}}
\def\bbr{{\mathbb R}}

\begin{document}

\title{Remarks on Finite Subset Spaces}

\author{Sadok Kallel}
\email{sadok.kallel@math.univ-lille1.fr}
\address{Laboratoire Painlev\'e\\
Universit\'e des Sciences et Technologies de Lille, France}

\author{Denis Sjerve}
\email{sjer@math.ubc.ca}
\address{Department of Mathematics\\
University of British Columbia, Vancouver, Canada}

\begin{abstract}
  This paper expands on and refines some known and less well-known
  results about the finite subset spaces of a simplicial complex $X$
  including their connectivity and manifold structure. It also
  discusses the inclusion of the singletons into the three fold subset
  space and shows that this subspace is weakly contractible but
  generally non-contractible unless $X$ is a cogroup.  Some
  homological calculations are provided.
\end{abstract}

\maketitle

\section{Statement of Results}\label{intro}

Let $X$ be a topological space (always assumed to be path connected),
and $k$ a positive integer. It has become increasingly useful in
recent years to study the space
$$\sub{n}X:= \{\{x_1,\ldots, x_\ell\}\subset X\ |\ \ell\leq n\}
$$ 
of all finite subsets of $X$ of cardinality at most $n$ \cite{akr,
  beilinson, handel, mp, rose, tuffley1}. This space is topologized as
the identification space obtained from $X^n$ by identifying two
$n$-tuples if and only if the sets of their coordinates coincide
\cite{borsuk}. The functors $\sub{n}(-)$ are homotopy functors in the
sense that if $X\simeq Y$ then $\sub{n}(X)\simeq\sub{n}(Y)$.  If $k\le
n$ then $Sub_kX$ naturally embeds in $Sub_nX.$ We write
$j_n:X\hookrightarrow Sub_nX$ for the inclusion given by
$j_n(x)=\{x\}$.

This paper takes advantage of the close relationship between finite
subset spaces and symmetric products to deduce a number of useful
results about them.

As a starting point, we discuss cell structures on finite subset
spaces. We observe in \S \ref{decomp} that if $X$ is a finite
$d$-dimensional simplicial complex, then $\sub{n}X$ is an
$nd$-dimensional CW complex and of which $\sub{k}X$ for $k\leq n$ is a
subcomplex (Proposition \ref{celldecomp}). Furthermore, $\sub{}X :=
\coprod_{n\geq 1}\sub{n}X$ has the structure of an abelian CW-monoid
(without unit) whenever $X$ is a simplicial complex.

In \S \ref{connect} we address a connectivity conjecture stated in
\cite{tuffley3}. We recall that a space $X$ is $r$-connected if
$\pi_i(X)=0$ for $i\leq r$. A contractible space is $r$-connected
for all positive $r$. In \cite{tuffley3} Tuffley proves that
$\sub{n}X$ is $n-2$ connected and conjectures that it is $n+r-2$
connected if $X$ is $r$-connected. We are able to confirm his
conjecture for the three fold subset spaces. In fact we show

\bth\label{main1} If $X$ is $r$-connected, $r\geq 1$ and $n\geq 3$,
then $\sub{n}X$ is $r+1$-connected.
\end{theorem}

In \S \ref{caveats} we address a somewhat surprising fact about the
embeddings $\sub{k}X\hookrightarrow\sub{n}X, k\leq n$. A theorem of
Handel \cite{handel} asserts that the inclusion $j : \sub{k}(X)
\hookrightarrow\sub{2k+1}(X)$ for any $k\geq 1$ is trivial on
homotopy groups (i.e. ``weakly trivial''). This is of course not
enough to conclude that $j$ is the trivial map, and in fact it need
not be. Let $\sub{k}(X,x_0)$ be the subspace of $\sub{k}X$ of all
finite subsets containing the basepoint $x_0\in X$. Handel's result
is deduced from the more basic fact that the inclusion $j_{x_0}:
\sub{k}(X,x_0)\hookrightarrow \sub{2k-1}(X,x_0)$ is weakly trivial.
The following theorem implies that these maps are often not
null-homotopic.

\bth\label{essential} The embeddings $X\hookrightarrow\sub{3}(X,x_0),
x\mapsto \{x,x_0\}$, and $j: X\hookrightarrow \sub{3}(X)$, $x\mapsto
\{x\}$, are both null-homotopic if $X$ is a cogroup.  If $X=S^1\times S^1$
is the torus, then both $j_3$ and $j_{x_0}$ are non-trivial in
homology and hence essential.
\end{theorem}

For a definition of a cogroup, see \S\ref{caveats}. In particular
suspensions are cogroups. The second half of Theorem \ref{essential}
follows from a general calculation given in \S\ref{caveats} which
exhibits a model for $\sub{3}(X,x_0)$ and uses it to show that its
homology is an explicit quotient of the homology of the symmetric
square $\sp{2}X$ by a submodule determined by the coproduct on
$H_*(X)$. One deduces in particular a homotopy equivalence between
$\sub{3}(\Sigma X,x_0)$ and the \textit{reduced} symmetric square
$\bsp{2}(\Sigma X)$ (cf. definition \ref{reducedcons} and proposition
\ref{equivalent}).  The methods in \S \ref{caveats} are taken up again
in \cite{sadok} where an explicit spectral sequence is devised to
compute $H_*(\sub{n}X)$ for any finite simplicial complex $X$ and any
$n\geq 1$.

The final section of this paper deals with manifold structures on
$\sub{n}X$ and top homology groups. It is known that
$\sub{2}X=\sp{2}X$ is a closed manifold if and only if $X$ is closed
of dimension $2$. This is a consequence of the fact that
$\sp{2}(\bbr^d)$ is not a manifold if $d>2$, while
$\sp{2}(\bbr^2)\cong\bbr^4$. The following complete description is
due to Wagner \cite{wagner}

\bth\label{main3} Let $X$ be a closed manifold of dimension $d\geq
1$. Then $\sub{n}X$ is a closed manifold if and only if either (i)
$d = 1$ and $n=3$, or (ii) $d=2$ and $n=2$.
\end{theorem}

This result is established in \S \ref{manifold} where we use in the
case $d\geq 2$ the connectivity result of theorem \ref{main1}, one
observation from \cite{mostovoy} and some homological calculations
from \cite{ks}. In the case $d=1$ we reproduce Wagner's cute
argument. Furthermore in that section we refine some results of
Handel \cite{handel} on the top homology groups of $\sub{n}X$ when
$X$ is a manifold. We point out that if $X$ is a closed orientable
manifold of dimension $d\geq 2$, then the top homology group
$H_{nd}(\sub{n}X)$ is trivial if $d$ is odd and is $\bbz$ if $d$ is
even. This group is always trivial if $X$ is not orientable
(see \S\ref{topdim}).

{\sc Acknowledgment}: This work was initiated at PIms in Vancouver
and the first author would like to thank the institute for its
hospitality.


\section{Basic Constructions}\label{basic}

All spaces $X$ in this paper are path connected, paracompact, and
have a chosen basepoint $x_0$.

The way we will think of $\sub{n}X$ is as a quotient of the $n$-th
symmetric product $\sp{n}X$. This symmetric product is the quotient
of $X^n$ by the permutation action of the symmetric group
$\mathfrak{S}_n$. The quotient map $\pi : X^n\lrar\sp{n}X$ sends
$(x_1,\ldots, x_n)$ to the equivalence class $[x_1,\ldots, x_n]$. It
will be useful sometimes to write such an equivalence class as an an
abelian product $x_1\ldots x_n$, $x_i\in X$. There are topological
embeddings \begin{equation}\label{basepoint} j_n: X\hookrightarrow
\sp{n}X\ \ \ ,\ \ \ x\mapsto xx_0^{n-1}
\end{equation}

The finite subset space $\sub{n}X$ is obtained from $\sp{n}X$
through the identifications
$$[x_1,\cdots , x_n]\ \sim\ [y_1, \cdots , y_n] \ \ \Longleftrightarrow\ \
\{x_1,\ldots, x_n\}=\{y_1,\ldots, y_n\}$$ In multiplicative
notation, elements of $\sub{n}X$ are products $x_1x_2\cdots x_k$
with $k\leq n,$  and subject to the identifications $x_1^2x_2\cdots
x_k \sim x_1x_2\cdots x_k$.

The topology of $\sub{n}X$ is the quotient topology inherited from
$\sp{n}X$ or $X^n$ \cite{handel}. When $X$ is Hausdorff this
topology is equivalent to the so-called {\sl Vietoris finite}
topology whose basis of open sets are sets of the form
$$[U_1,\ldots, U_k]:=\{S\in\sub{n}X\ |\
S\subset\bigcup_{i=1}^kU_i\ \hbox{and}\ S\cap U_i\neq\emptyset\
\hbox{for each $i$}\}$$ where $U_i$ is open in $X$ \cite{wagner}.
When $X$ is a metric space,  $\sub{k}X$ is again a metric space
under the Hausdorff metric, and hence inherits a third and
equivalent topology \cite{wagner}. In all cases, for any topology we
use,  continuous maps between spaces induce continuous maps between
their finite subset spaces.

\bex\label{e3s2} Of course $\sub{1}X = X$ and $\sub{2}X=\sp{2}X$.
Generally, if ${\bf\Delta}^{n+1}X \subset\sp{n+1}X$ denotes the
image of the fat diagonal in $X^{n+1}$; that is
$${\bf\Delta}^{n+1}X :=\{x_1^{i_1}\ldots x_r^{i_r}\in\sp{n+1}X\ |\
r\leq n, \sum i_j = n+1\ 
\hbox{and}\ \ i_j>0\}$$ then there is a map $q : {\bf\Delta}^{n+1}X
\lrar\sub{n}X$, $x_1^{i_1}\ldots x_r^{i_r}\lrar \{x_1,\ldots ,x_r\}$, and
a pushout diagram
\begin{equation}\label{pushout2}
\xymatrix{
{\bf\Delta}^{n+1}X\ar[r]^{i}\ar[d]^q&\sp{n+1}X\ar[d]\\
\sub{n}X\ar[r]&\sub{n+1}X
}
\end{equation}
This is quite clear since we obtain $\sub{n+1}X$ by identifying
points in the fat diagonal to points in $\sub{n}X$. In particular,
when $n=2$, we have the pushout
\begin{equation}\label{pushout3}
\xymatrix{
X\times X\ar[r]^{i}\ar[d]^q&\sp{3}X\ar[d]\\
\sp{2}X\ar[r]&\sub{3}X
}
\end{equation}
where $q(x,y)=xy$ and $i(x,y)= x^2y$. The homology of $\sub{3}(X)$ can
then be obtained from a Mayer-Vietoris sequence. Some calculations for
the three fold subset spaces are in \S \ref{caveats}.

There are two immediate and non-trivial consequences of the above
pushouts.  Albrecht Dold shows in \cite{dold} that the homology of the
symmetric products of a CW complex $X$ only depends on the homology of
$X$. The pushout diagram in (\ref{pushout2}) shows that in the case of
the finite subset spaces, this homology also depends on the
\textit{cohomology structure of $X$}. This general fact for the three
and four fold subset spaces is further discussed in \cite{taamallah}.

The second consequence of (\ref{pushout2}) is that it yields an
important corollary.\eex

\bco\label{fundamental} $\sub{n}X$ is simply connected for $n\geq
3$. \eco

\begin{proof}
We use the following known facts about symmetric products:
$\pi_1(\sp{n}X)\cong H_1(X; \bbz )$ whenever $n\geq 2$, and the
inclusion $j_n: X\hookrightarrow\sp{n}X$ induces the abelianization
map at the level of fundamental groups (P.A. Smith \cite{smith} proves
this for $n=2,$ but his argument applies for $n>2$
\cite{taamallah}). For $n\geq 3$, consider the composite
$$X\fract{\alpha}{\lrar}
{\bf\Delta}^{n}X\fract{i}{\lrar}\sp{n}X$$ with $\alpha (x) =
[x,x_0,\ldots, x_0]$. The induced map $j_{n*}=i_*\circ\alpha_*$ on
$\pi_1$ is surjective, as we pointed out, and hence so is $i_*$.
Assume we know that $\pi_1(\sub{3}(X))=0$. Then the fact that $i_*$
is surjective implies immediately by the Van-Kampen theorem and the
pushout diagram in (\ref{pushout2}) that $\pi_1(\sub{4}X)=0.$
  By induction we see that $\pi_1(\sub{n}X)=0$ for larger $n$.
Therefore, we need only establish the claim for $n=3.$  For that we
apply Van Kampen to diagram (\ref{pushout3}).  Consider the maps
$\tau : x_0\times X \hookrightarrow X\times
X\fract{i}{\lrar}\sp{3}X$ and $\beta : X\times x_0\lrar X\times
X\fract{q}{\lrar}\sp{2}X.$ Now $i(x,y)=x^2y$ so that $\tau (x_0,x) =
x_0^2x = j_3(x)$ and $\beta (x,x_0) = xx_0 = j_2(x).$ Since the
$j_k$'s are surjective on $\pi_1$ it follows that $\tau$ and $\beta$
are surjective on $\pi_1.$  Therefore, for any classes
$u\in\pi_1(SP^3X)$ and $v\in\pi_1(SP^2X)$, $\exists$ a class
$w\in\pi_1(X\times X)$ such that $i_*(w)=u$ and $q_*(w)=v.$ This
shows that $\pi_1(Sub_3X)=0.$
\end{proof}

This corollary also follows from \cite{biro, tuffley3}, where it is
shown that $\sub{n}X$ is $(n-2)$-connected for $n\geq 3$. However, the
proof above is completely elementary.

\subsection{Reduced Constructions}\label{reducedcons} For the
spaces under consideration, the natural inclusion
$\sub{n-1}X\subset\sub{n}X$ is a cofibration \cite{handel}. We write
$\bsub{n}X:=\sub{n}X/\sub{n-1}X$ for the cofiber. Similarly
$\sp{n-1}X$ embeds in $\sp{n}X$ as the closed subset of all
configurations $[x_1,\ldots , x_n]$ with $x_i$ at the basepoint for
some $i$. We set $\bsp{n}X := \sp{n}X/\sp{n-1}X$.

Note that even though $\sp{2}X$ and $\sub{2}X$ are the same, there is an
essential difference between their reduced analogs. The difference
here comes from the fact that the inclusion
$X\hookrightarrow\sub{2}X$ is the composite $X\fract{\Delta}{\lrar}
X\times X\lrar\sp{2}X\cong \sub{2}X,$ where $\Delta$ is the
diagonal, while $j_2: X\hookrightarrow\sp{2}X$ is the basepoint
inclusion.

\bex\label{connectbsp2s} When $X=S^1$, $\sp{2}(S^1)$ is the closed
M\"{o}bius band. If we view this band as a square with two sides
identified along opposite orientations, then
$S^1=\sp{1}(S^1)\hookrightarrow\sp{2}(S^1)$ embeds into this band as
an edge (see figures on p. 1124 of \cite{tuffley1}).  Hence this
embedding is homotopic to the embedding of an equator, and so
$\bsp{2}(S^1)$ is contractible. On the other hand $S^1=\sub{1}(S^1)$
embeds into $\sub{2}(S^1)=\sp{2}(S^1)$ as the diagonal $x\mapsto
\{x,x\}=[x,x],$ which is the boundary of the M\"{o}bius band, and so
$\bsub{2}(S^1)=\bbr P^2$. \eex

\bex When $X=S^2$, $\sp{2}(S^2)$ is the complex projective plane
$\bbp^2,$
  $\sp{1}(S^2)=\bbp^1$ is a hyperplane, and
$\bsp{2}(S^2)=S^4$. On the other hand $\bsub{2}(S^2)$ has the
following description. Write $\bbp^1$ for $\bbc\cup\{\infty\}$. Then
$\bsub{2}(S^2)$ is the quotient of $\bbp^2$ by the image of the
Veronese embedding $\bbp^1\lrar\bbp^2$, $z\mapsto [z^2 : -2z :1]$,
$\infty\mapsto [1:0:0]$. To see this, identify $\sp{n}(\bbc )$ with
$\bbc^n$ by sending $(z_1,\ldots, z_n)$ to the coefficients of the
polynomial $(x-z_1)\ldots (x-z_n)$. This extends to the
compactifications to give an identification of $\sp{n}(S^2)$ with
$\bbp^n$ (\cite{hatcher}, chapter 4). When $n=1$, $(z,z)$ is mapped
to the coefficients of $(x-z)(x-z),$ that is to $(z^2,-2z)$. Note
that the diagonal $S^2\lrar\sp{2}(S^2)=\bbp^2$ is multiplication by
$2$ on the level of $H_2$ so that, in particular,
$H_4(\bsub{2}(S^2))=\bbz,$  $H_2(\bsub{2}(S^2))=\bbz_2,$ and all
other reduced homology groups are zero. \eex


\section{Cell Decomposition}\label{decomp}

If $X$ is a simplicial complex, there is a standard way to pick a
$\mathfrak S_n$-equivariant simplicial decomposition for the product
$X^n$ so that the quotient map $X^n\lrar\sp{n}X$ induces a cellular
structure on $\sp{n}X$. We argue that this same cellular structure
descends to a cell structure on $\sub{n}X$. The construction of this
cell structure for the symmetric products is fairly classical
\cite{liao, nakaoka}. The following is a review and slight
expansion.

\bpr\label{celldecomp} Let $X$ be a simplicial complex.  For $n\geq
1$ there exist cellular decompositions for $X^n$, $\sp{n}X$ and
$\sub{n}X$ so that all of the quotient maps
$X^n\rightarrow\sp{n}X\rightarrow\sub{n}X$ and the concatenation
pairings $+$ are cellular
\begin{equation}\label{pairing}
\xymatrix{
\sp{r}X\times\sp{s}X\ar[r]^{\ \ +}\ar[d]&\sp{r+s}X\ar[d]\\
\sub{r}X\times\sub{s}X\ar[r]^{\ \ +}&\sub{r+s}X
}
\end{equation}
Furthermore the subspaces ${\bf\Delta}^{n},\sp{n-1}X\subset \sp{n}X$
and $\sub{n-1}X \subset\sub{n}X$ are subcomplexes.
\epr

\begin{proof}
Both $\sp{n}X$ and $\sub{n}X$ are obtained from $X^n$ via
identifications. If for some simplicial (hence cellular) structure
on $X^n$, derived from that on $X$, these identifications become
simplicial (i.e. they identify simplices to simplices), then the
quotients will have a cellular structure and the corresponding
quotient maps will be cellular with respect to these structures.

As we know, one obtains a nice and natural $\sn$-equivariant
simplicial structure on the product if one works with \textit{ordered}
simplicial complexes \cite{liao, nakaoka, dwyer}. We write $X_\bullet$
for the abstract simplicial (i.e. triangulated) complex of which $X$
is the realization. So we assume $X_\bullet$ to be endowed with a
partial ordering on its vertices which restricts to a total ordering
on each simplex. Let $\prec$ be that ordering. A point $w =
(v_1,\ldots, v_n)$ is a vertex in $X^n_\bullet$ if and only if $v_i$
is a vertex of $X_\bullet$. Different vertices
\begin{equation}\label{vertices}
w_0 = (v_{01}, v_{02}, \ldots, v_{0n})\ ,\ldots,\  w_k = (v_{k1},
v_{k2},\ldots, v_{kn})
\end{equation}
span a $k$-simplex in $X^n_\bullet$ if, and only if, for each $i$,
the $k+1$ vertices $v_{0i},v_{1i},\ldots, v_{ki}$ are contained in a
simplex of $X$ and $v_{0i}\prec v_{1i}\prec\cdots\prec v_{ki}$. We
write  $\varpi:=[w_0, \ldots, w_k ]$ for such a simplex.

The permutation action of $\tau\in\sn$ on $\varpi=[w_0, \ldots, w_k
]$ is given by $\tau \varpi = [\tau w_0, \ldots, \tau w_n ]$. This
is a well-defined simplex since the factors of each vertex $w_j =
(v_{j_11}, v_{j_22},\ldots, v_{j_nn})$ are permuted simultaneously
according to $\tau,$ and hence the order $\prec$ is preserved.  The
permutation action is then simplicial and $\sp{n}X$ inherits a CW
structure by passing to the quotient.

{\bf Fact 1}: If a point $p:= (x_1,x_2,\ldots, x_{n})\in X^n$ is such
that $x_{i_1}=x_{i_2}=\ldots = x_{i_r}$, then $p$ lies in some
$k$-simplex $\varpi$ whose vertices $[w_0,\ldots, w_k]$ are such that
$v_{ji_1}=v_{ji_2}=\cdots =v_{ji_r}$ for $j=0,\ldots, k$. This implies
that the fat diagonal is a simplicial subcomplex. It also implies that
any permutation that fixes such a point $p$ must fix the vertices of
the simplex it lies in and hences fixes it pointwise. In other words, if a
permutation leaves a simplex invariant then it must fix it pointwise.

{\bf Fact2}: If $p=(x_1,x_2,\ldots, x_n)\in\varpi$ is a simplex with
vertices $w_0, ..., w_k$ as in (\ref{vertices}), and if $\pi :
X^n\lrar X^i$ is any projection, then $\pi (p)$ lies in the simplex
with vertices $\pi (w_0),\cdots, \pi (w_k)$ (which may or may not be
equal). For instance $\pi (p):=(x_1,\ldots, x_i)$ lies in the simplex
with vertices $(v_{01}, v_{02},\ldots, v_{0i})$, ..., $(v_{k1},
v_{k2},\ldots, v_{ki})$.

We are now in a position to see that $\sub{n}X$ is a CW complex.
Recall that $\sub{n}X=X^n/_\sim$ where
$$(x_1,\ldots, x_n)\sim (y_1,\ldots, y_n)\Longleftrightarrow
\{x_1,\ldots, x_n\} = \{y_1,\ldots, y_n\}$$ Clearly, if
$(x_1,\ldots, x_n)\sim (y_1,\ldots, y_n)$ then $\tau (x_1,\ldots,
x_n)\sim \tau (y_1,\ldots, y_n)$ for $\tau\in\mathfrak S_n$. We wish
to show that these identifications are simplicial. Let's argue
through an example (the general case being identical). We have the
identifications in $Sub_6X$:
\begin{equation}\label{identify}
p:=(x,x,x,y,y,z)\sim (x,x,y,y,y,z)=:q
\end{equation}
By using Fact 2 applied to the projection skipping the third
coordinate and then Fact 1, we can see that $p$ and $q$ lie in
simplices with vertices of the form $(v_1,v_1,?,v_2,v_2,v_3)$. By
using Fact 1 again, $p$ lies in a simplex $\sigma_p$ with vertices of
the form $(v_1,v_1,v_1,v_2,v_2,v_3)$ while $q$ lies in a simplex
$\sigma_q$ with vertices of the form $(v_1,v_1,v_2,v_2,v_2,v_3)$. It
follows that the identification (\ref{identify}) identifies vertices
of $\sigma_p$ with vertices of $\sigma_q,$ and hence identifies
$\sigma_p$ with $\sigma_q$ as desired.

In conclusion, the quotient $\sub{n}X$ inherits a cellular structure
and the composite
$$X^n\fract{\pi}{\lrar}\sp{n}X\fract{q}{\lrar}\sub{n}X$$
is cellular. Since the pairing (\ref{pairing}) is covered by
$X^r\times X^s\lrar X^{r+s},$ which is simplicial (by construction),
and since the projections are cellular, the pairing (\ref{pairing})
must be cellular.
\end{proof}

\bre We could have worked with simplicial sets instead \cite{biro}.
Similarly, Mostovoy (private communication) indicates how to
construct a simplicial set $\sub{n}X$ out of a simplicial set $X$
such that $|\sub{n}X| = \sub{n}|X|$. This approach will be further
discussed in \cite{sadok}.\ere

  The following corollary is also obtained in \cite{biro}.

\bco\label{topcell} For $X$ a simplicial complex, $\sub{k}X$ has a
CW decomposition with top cells in $k\dim X,$ so that $H_*(\sub{k}X)
= 0$ for $*> k\dim X$.  \eco

We collect a couple more corollaries

\bco\label{same} If $X$ is a $d$-dimensional complex with $d\geq 2$,
then the quotient map $\sp{n}X\to\sub{n}X$ induces a homology
isomorphism in top dimension $nd$. \eco

\begin{proof} When $X$ is as in the hypothesis,
$\sub{n-1}X$ is a codimension $d$ subcomplex of $\sub{n}X$ and since
$d\geq 2$, $H_{nd}(\sub{n}X) = H_{nd}(\sub{n}X,\sub{n-1}X)$. On the
other hand, Proposition \ref{celldecomp} implies that
${\bf\Delta}^{n}X $ is a codimension $d$ subcomplex of $\sp{n}X$ so
that $H_{nd}(\sp{n}X)\cong H_{nd}(\sp{n}X,{\bf\Delta}^{n}X )$ as
well. But according to diagram (\ref{pushout2}), we have the
homeomorphism
$$\sp{n}X/{\bf\Delta}^{n}X \cong \sub{n}X/\sub{n-1}X
$$
Combining these facts yields the claim.
\end{proof}

\bco\label{sameconnect} Both $\sp{k}X$ and the fat diagonal
${\bf\Delta}^k\subset \sp{k}X$ have the same connectivity as $X$,
and this is sharp. \eco

\begin{proof} If $X$ is an $r$-connected ordered simplicial
complex, then $X$ admits a simplicial structure so that the
$r$-skeleton $X_r$ is contractible in $X$ to some point $x_0\in X$.
With such a simplicial decomposition we can consider Liao's induced
decomposition $X^k_{\bullet}$ on $X^k$ and its $r$-skeleton $X^k_r$.
Note that
$$X^k_r \subset \bigcup_{i_1+\cdots + i_k\leq r}
X_{i_1}\times X_{i_2}\times\cdots\times X_{i_k} \subset (X_r)^k$$ If
$F :X_r\times I\lrar X$ is a deformation of $X_r$ to $x_0$, then
$F^k$ is a deformation of $(X_r)^k$, hence $X^k_r$, to $(x_0,\ldots,
x_0)$ in $X^k,$ and this deformation is $\mathfrak{S}_k$
equivariant. Since the $r$-skeleton of $\sp{k}X$ is the
$\mathfrak{S}_k$-quotient of $X^k_r$, it is then itself contractible
in $\sp{k}X,$ and this proves the first claim. Similarly, the
simplicial decomposition we have introduced on $X^k$ includes the
fat diagonal $\Lambda^k$ as a subcomplex with $r$-skeleton
$\Lambda_r^k := \Lambda^k\cap X^k_r$. The deformation $F^k$
preserves the fat diagonal  and so it restricts to $\Lambda^k$ and
to an equivariant deformation $F^k : \Lambda^k_r\times
I\lrar\Lambda^k$. This means that the $r$-skeleton of
$q(\Lambda^k)=:{\bf\Delta}^k\subset\sp{k}X$ is itself contractible
in ${\bf\Delta}^k,$ and the second claim follows. This bound is
sharp for symmetric products since when $X=S^2$,
$\sp{2}(S^2)=\bbp^2$. It is sharp for the fat diagonal as well since
${\bf\Delta}^3X\cong X\times X$ has exactly the same connectivity of
$X$.
\end{proof}


\section{Connectivity}\label{connect}

As we've established in corollary  \ref{fundamental}, finite subset
spaces $Sub_nX,\ n\ge 3,$ are always simply connected. In this
section we further relate the connectivity of $\sub{k}X$ to that of
$X$. We first need the following useful result proved in
\cite{braid}.

\bth\label{connectbsp} If $X$ is $r$-connected with $r\geq 1$, then
$\bsp{n}X$ is $2n+r-2$ connected.
\end{theorem}

Example \ref{bsp2} shows that $\bsp{2}(S^k)$ is $k+1$-connected as
asserted. Note that $\bsp{2}(S^2)=S^4$ is $3$-connected, so theorem
\ref{connectbsp} is sharp.

\bco\label{nakak} (\cite{nakaoka} corollary 4.7) If $X$ is
$r$-connected, $r\ge 1,$ then $H_{*}(X)\cong H_{*}(\sp{n}X)$ for
$*\leq r+2$. This isomorphism is induced by the map $j_n$ adjoining
the base point.\eco

\begin{proof} We give a short proof based on theorem
\ref{connectbsp}. By Steenrod's homological splitting \cite{nakaoka}
\begin{equation}\label{steenrod}
H_*(\sp{n}X)\cong \bigoplus_{k=1}^n H_*(\sp{k}X,\sp{k-1}X) =
\bigoplus_{k=2}^n \tilde H_*(\bsp{k}X)\oplus H_*(X)
\end{equation} with
$\sp{0}X=\emptyset$. But $\tilde H_*(\bsp{k}X)=0$ for $*\leq
2k+r-2$. The result follows.
\end{proof}

\bre Note that corollary \ref{nakak} cannot be improved to $r=0$
(i.e. $X$ connected). It fails already for the wedge $X=S^1\vee S^1$
and $n=2$ since $\sp{2}(S^1\vee S^1)\simeq S^1\times S^1$ (see
\cite{ks}) and hence $H_2(\sp{2}(S^1\vee S^1))\not\cong H_2(S^1\vee S^1)$.
Note also that (\ref{steenrod}) implies that $H_*(X)$ embeds into
$H_*(\sp{n}X)$ for all $n\geq 1$; a fact we will find useful
below.\ere

\bpr\label{connectha} Suppose $X$ is $r$-connected, $r\geq 1$. Then
$\sub{k}X$ is $r+1$ connected whenever $k\geq 3$.  \epr

\begin{proof} Write $x_0\in X$ for the basepoint and assume $k\geq 3$.
Remember that the $\sub{k}X$ are simply connected for $k\geq 3$
(corollary \ref{fundamental}) so by the Hurewicz theorem if they
have trivial homology up to degree $r+1$, then they are connected up
to that level. We will now show by induction that $H_*(\sub{k}X)=0$
for $*\leq r+1$. The first step is to show that
$H_*(\sp{k}X,{\bf\Delta}^k)=H_*(\sub{k}X,\sub{k-1}X)=0$ for $* \leq
r+1$. We write $i: {\bf\Delta}^k\hookrightarrow\sp{k}X$ for the
inclusion.

From the fact that   ${\bf\Delta}^k$ and $\sp{k}X$ have the same
connectivity as $X$ (corollary \ref{sameconnect}), their homology
vanishes up to degree $r$ which implies similarly that the relative
groups are trivial up to that degree. On the other hand $X$ embeds
in ${\bf\Delta}^k$ via $x\mapsto [x,x_0,\cdots, x_0]$ (this is a
well-defined map since $k\geq 3$) and, since the composite $j_k :
X\rightarrow {\bf\Delta^k}\fract{i}{\lrar} \sp{k}X$ is an
isomorphism on $H_{r+1}$ (corollary \ref{nakak}), we see that the
map $i_*: H_{r+1}({\bf\Delta}^k)\lrar H_{r+1}(\sp{k}X)$ is
surjective. Hence $H_{r+1}(\sp{k}X, {\bf\Delta}^k)= 0$.

Now since $0=H_*(\sp{k}X,{\bf\Delta}^k)= H_*(\sub{k}X,\sub{k-1}X)$
for $*\leq r+1$, it follows that $H_*(\sub{k-1}X)\cong
H_*(\sub{k}X)$ for $*\leq r$ and that $H_{r+1}(\sub{k-1}X)\lrar
H_{r+1}(\sub{k}X)$ is surjective. So if we prove that
$H_*(\sub{3}X)=0$ for $*\leq r+1$, then by induction we will have
proved our claim.

Consider the homology long exact sequences for $(\sub{3}X,\sub{2}X)$
and $(\sp{3}X,{\bf\Delta}^3X),$ where again we identify
${\bf\Delta}^3X$ with $X\times X$.  We obtain commutative diagrams
$$\xymatrix{
\ar[r]&H_{r+2}(\sub{3}X,\sub{2}X)\ar[r]&
H_{r+1}(\sub{2}X)\ar[r]^{i_*}&H_{r+1}(\sub{3}X)\ar[r]&0\\
\ar[r]&H_{r+2}(\sp{3}X,X^2)\ar[r]\ar[u]^\cong&
H_{r+1}(X^2)\ar[r]^{\alpha_*}\ar[u]^{q_*}&H_{r+1}(\sp{3}X)\ar[r]\ar[u]^{\pi_*}&0
}
$$
where $\alpha (x,y)=x^2y$ and $\pi : \sp{3}X\lrar\sub{3}X$ is the
quotient map.  We want to show that $i_* = 0$ so that by exactness
$H_{r+1}(\sub{3}X)=0$.  Now $q_*$ is surjective since the composite
$$X\lrar X\times \{x_0\}
\hookrightarrow X\times X\lrar \sp{2}X=\sub{2}X$$ induces an
isomorphism on $H_{r+1}$ by Corollary \ref{nakak}.  Showing that
$i_*=0$ comes down therefore to showing that $\pi_*\circ\alpha_*
=0$. But note that for $r\geq 1,$ which is the connectivity of $X$,
classes in $H_{r+1}(X\times X)$ are necessarily spherical and we
have the following commutative diagram
$$\xymatrix{
\pi_{r+1}X\times\pi_{r+1}(X)\ar[r]^{\ \ \cong}&\pi_{r+1}(X\times X)\ar[rr]\ar[d]^h&&
\pi_{r+1}(\sub{3}(X))\ar[d]^{h}\\
&H_{r+1}(X\times X)\ar[rr]^{\pi_*\circ\alpha_*}&&H_{r+1}(\sub{3}(X))
}$$
where $h$ is the Hurewicz homomorphism. The top map is trivial since
when restricted to each factor $\pi_{r+1}(X)$ it is trivial
according to the useful theorem \ref{handel} below (or to
corollary \ref{cowt}) . Since $h$ is surjective,
$\pi_*\circ\alpha_*=0$ and $H_{r+1}(\sub{3}X) = 0$ as desired.
\end{proof}


\section{The Three Fold Finite Subset Space}\label{caveats}

There are many subtle points that come up in the study of finite
subset spaces. We illustrate several of them through the study of
the pair $(\sub{3}X, X)$. The three fold subset space has been
studied in \cite{mostovoy, rose, tuffley1} for the case of the
circle and in \cite{tuffley2} for topological surfaces.

Again all spaces below are assumed to be connected. We say a map is
weakly contractible (or weakly trivial) if it induces the trivial map
on all homotopy groups. The following is based on a cute argument well
explained in \cite{handel} or (\cite{beilinson} section 3.4).

\bth\label{handel}\cite{handel} $\sub{k}(X)$ is weakly contractible
in $\sub{2k+1}(X)$.
\end{theorem}

{\sc Caveat 1}: A map $f : A\lrar Y$ being weakly contractible does
not generally imply that $f$ is null homotopic. Indeed let $T$ be
the torus and consider the projection $T\lrar S^2$ which collapses
the one-skeleton. Then this map induces an isomorphism on $H_2$ but
is trivial on homotopy groups since $T=K(\bbz^2,1)$. Of course if
$A=S^k$ is a sphere, then ``weakly trivial'' and ``null-homotopic''
are the same since the map $A\lrar Y$ represents
  the zero element in $\pi_{k}Y$. For example, in (\cite{ch},
lemma 3.3), the authors construct explicitly an extension of the
inclusion $S^n\hookrightarrow\sub{3}(S^n)$ to the disk
$B^{n+1}\lrar\sub{3}(S^n)$, $\partial B^{n+1}=S^n$. This section
argues that this implication doesn't generally hold for
non-suspensions.

\vskip 7pt {\sc Caveat 2}: When comparing symmetric products to
finite subset spaces, one has to watch out for the fact that the
basepoint inclusion $\sp{k}(X)\lrar\sp{k+1}(X)$ {\it does not
commute} via the projection maps with the inclusion
$\sub{k}(X)\lrar\sub{k+1}(X)$. This has already been pointed out in
example \ref{connectbsp2s} and is further illustrated in the
corollary below.

\bco\label{cowt} 
The composite $\sp{k}(X)\lrar\sp{2k+1}(X)\lrar\sub{2k+1}(X)$ is
weakly trivial. \eco

\begin{proof} This map is equivalent to the composite
\begin{equation}\label{composite}
\sp{k}(X)\lrar\sub{k}(X)\fract{\mu
}{\lrar}\sub{k+1}(X,x_0)\hookrightarrow\sub{2k+1}(X)
\end{equation}
where $\mu (\{x_1,\ldots, x_k\}) =\{x_0,x_1,\ldots, x_k\}$, $x_0$ is
the basepoint of $X$ and $\sub{k+1}(X,x_0)$ is the subspace of
$\sub{k+1}(X)$ of all subsets containing this basepoint. Note that
$\mu$ is not an embedding as pointed out in \cite{tuffley2} but is
one-to-one away from the fat diagonal. The key point here is
again (\cite{handel}, Theorem 4.1) which asserts that the inclusion
$$\sub{k+1}(X,x_0)\hookrightarrow\sub{2k+1}(X,x_0)$$
is weakly contractible. This in turn implies that the last map in
(\ref{composite}) is weakly trivial as well and the claim follows.
\end{proof}

{\sc Caveat 3}: For $n\geq 2$, one can embed
$X\hookrightarrow\sub{n}(X)$ in several ways. There is of course the
natural inclusion $j$ giving $X$ as the subpace of singletons. There
is also, for any choice of $x_0\in X,$ the embedding $j_{x_0}:
x\mapsto \{x,x_0\}$. Any two such embeddings for different choices
of $x_0$ are equivalent when $X$ is path-connected (any choice of a
path between $x_0$ and $x'_0$ gives a homotopy between $j_{x_0}$ and
$j_{x'_0}$). It turns out however that $j$ and $j_{x_0}$ are
fundamentally different. The simplest example was already pointed
out for $S^1,$ where $\sub{2}(S^1)$ was the M\"{o}bius band with $j$
being the embedding of the boundary circle while $j_{x_0}$ is the
embedding of an equator.

One might ask the question whether it is true that $j$ is
null-homotopic if and only if $j_{x_0}$ is null-homotopic? This is
at least true for suspensions as the next lemma illustrates.

Recall that a co-$H$ space $X$ is a space whose diagonal map factors
up to homotopy through the wedge; that is there exists a $\delta$
such that the composite
$$X\fract{\delta}{\lrar} X\vee X\hookrightarrow X\times X$$
is homotopic to the diagonal $\Delta : X\lrar X\times X, x\mapsto
(x,x)$. A cogroup $X$ is a co-$H$ space that is co-associative with
a homotopy inverse. This latter condition means there is a map $c:
X\lrar X$ such that $X\fract{\delta}{\lrar}X\vee X\fract{c\vee
1}{\lrar} X$ is null-homotopic. This is in fact the definition of a
left inverse but it implies the existence of a right inverse as well
\cite{arkowitz}. If $X$ is a cogroup, then for every based space
$Y$, the set of based homotopy classes of based maps $[X,Y]$ is a
group. The suspension of a space is a cogroup and there exist
several interesting cogroups that are not suspensions
(\cite{arkowitz}, \S4).

Write $j_{x_0}: X\hookrightarrow\sub{3}(X,x_0)$ the map
$x\mapsto\{x,x_0\}$. Its continuation to $\sub{3}(X)$ is also
written $j_{x_0}$.

\ble\label{null} Suppose $X$ is a cogroup. Then the embeddings $j_{x_0}:
X\hookrightarrow\sub{3}(X,x_0)$ and $j: X\hookrightarrow\sub{3}(X)$
are null-homotopic. \ele

\begin{proof} The argument in \cite{handel} extends to this situation.
We deal with $j_{x_0}$ first. This is a based map at $x_0.$ Its
homotopy class $[j_{x_0}]$ lives in the group
$G=[X,\sub{3}(X,x_0)]$. The following composite is checked to be
again $j_{x_0}$.
$$j_{x_0} : X\fract{\Delta}{\lrar}X\times
X\fract{j_{x_0}+j_{x_0}}{\ra 4}\sub{3}(X,x_0)$$ This factors up to
homotopy through the wedge $\iota : X\fract{\delta}{\lrar}X\vee
X\fract{j_{x_0}\vee j_{x_0}}{\ra 3}\sub{3}(X,x_0)$. Of course
$[\iota ]=[j_{x_0}]$. But observe that $[\iota]=2[j_{x_0}]$ by
definition of the additive structure of $G$. This means that
$[j_{x_0}]=2[j_{x_0}];$ thus $[j_{x_0}]=0$ and $j_{x_0}$ is trivial
(through a homotopy fixing $x_0$)

Let's now apply this to the inclusion $j: X\hookrightarrow\sub{3}(X)$
which is assumed to be based at $x_0$. We also denote the composite
$X\fract{j_{x_0}}{\lrar}\sub{3}(X,x_0)\lrar\sub{3}X$ by $j_{x_0}$. Using the
co-H structure as before we get the commutative diagram
$$\xymatrix{
X\ar[rr]^{\Delta}\ar[d]^\delta&&X\times X\ar[d]^{j+j}\\
X\vee X\ar[rr]^{j_{x_0}\vee j_{x_0}}&& \sub{3}(X)
}
$$
Since $j_{x_0}$ was just shown to be null homotopic, then so is $j =
(j+ j)\circ \Delta$.
\end{proof}

Let's now turn to the second part of theorem \ref{essential}.

\subsection{The Space $\sub{3}(X,x_0)$}
The preceeding discussion shows the usefulness of looking at the based
finite subset space $\sub{n}(X,x_0)$. We start with a key
computation. Write $\Delta$ for the diagonal $X\lrar\sp{2}X$,
$x\mapsto [x,x]$, and identify the image of $j_*:
H_*(X)\hookrightarrow H_*(\sp{2}(X))$ with $H_*(X)$ by the Steenrod
homological splitting (\ref{steenrod}).

\ble\label{computation} Let $X$ be a compact cell complex. Then
$H_*(\sub{3}(X,x_0)) = H_*(\sp{2}X)/I$ where $I$ is the submodule
generated by $\Delta_*c-c, c\in H_*(X)\hookrightarrow H_*(\sp{2}X)$.
\ele

\begin{proof}
Start with the map $\alpha : \sp{2}(X)\lrar\sub{3}(X,x_0),
[x,y]\mapsto \{x,y,x_0\}$ which is surjective and generically
one-to-one (i.e. one-to-one on the subspace of points $[x,y]$ with
$x\neq y$). Observe that $\alpha ([x,x]) = \alpha ([x,x_0])$. This
implies that $\sub{3}(X,x_0)$ is homeomorphic to the identification
space
\begin{equation}\label{wahid}
\sp{2}(X)/\sim \ \ , \ \ [x,x]\sim [x,x_0]\ \ , \ \forall x\in X
\end{equation} In order to compute the homology of this quotient we
will replace it with the following space
\begin{equation}\label{thneen}
W_2(X) := \sp{2}(X)\sqcup X\times I/\sim\ \ ,\ \ [x,x]\sim (x,1)\ ,\
[x,x_0]\sim (x,0)\ ,\ [x_0,x_0]\sim (x_0,t)
\end{equation}
It is not hard to see that (\ref{wahid}) and (\ref{thneen}) are
homotopy equivalent. We can easily see that these spaces are homology
equivalent as follows (this is enough for our purpose).  There is a
well-defined map $g : W_2(X)\lrar\sp{2}(X)/\sim$ sending $[x,y]\mapsto
[x,y], (x,t)\mapsto [x,x_0]$. The inverse image $g^{-1}([x,y])=[x,y]$
if $x\neq y$ and both points are different from $x_0$. The inverse
image of $[x,x]$ or $[x,x_0]$ is an interval when $x\neq x_0$, hence
contractible, and it is a point when $x=x_0$. In all cases preimages
under $g$ are acyclic and hence $g$ is a homology equivalence by the
Begle-Vietoris theorem. The homology structure of $\sub{3}(X,x_0)$ can
be made much more apparent using the form (\ref{thneen}) and this is
why we have introduced it.

Let $(C_*(\sp{2}(X)),\partial) $ be a chain complex for $\sp{2}(X)$
containing $C_*(X)$ as a subcomplex and for which the diagonal map
$X\lrar\sp{2}X$ is cellular. Associate to $c\in C_i(X)$ a chain $|c|$
in degree $i+1$ representing $I\times c\in C_{i+1}(I\times X)$ if
$c\neq x_0$ (the $0$-chain representing the basepoint). We write
$|C_*(X)|$ for the set of all such chains. The geometry of our
construction gives a chain complex for $W_2(X)$ as follows
\begin{equation}\label{complex}
C_*(W_2(X)) = C_*(\sp{2}(X))\oplus |C_*(X)|
\end{equation}
with boundary $d$ such that $d(c)=\partial c$ and
$$d|c| = c - \Delta_*(c)-|\partial c| $$ This comes from the formula
for the boundary of the product of two cells which is in general given
by $\partial (\sigma_1\times\sigma_2) =\partial
(\sigma_1)\times\sigma_2 + (-1)^{|\sigma_1|}\sigma_1\times \partial
(\sigma_2)$. We check indeed that $d\circ d=0$.  To compute the
homology we need to understand cycles and boundaries in this chain
complex. Write a general element of (\ref{complex}) as $\alpha+
|c|$. The boundary of this element is $\partial \alpha + c -
\Delta_*(c)-|\partial c|, $ and it is zero if, and only if, $\partial
\alpha = \Delta_*(c) - c$ and $|\partial c|=0$. That is if, and only
if, $c$ is a cycle and $\Delta_*(c) - c$ is a boundary.  This means
that in $H_*(\sp{2}(C))$, $\Delta_*(c) = c$. We claim this is not
possible unless $c=0$. Indeed, if $c$ is a positive dimensional
(homology) class, then $\Delta_*(c) = c\tensor 1 + \sum c'\tensor c''
+ 1\tensor c$ in $H_*(X\times X)$ and hence in $H_*(\sp{2}(C))$,
$\Delta_*(c) = 2c + \sum c'*c''$ where by definition $c'*c'' =
q_*(c'\tensor c'')$, $q : X\times X\lrar\sp{2}(X)$ the
projection. This can never be equal to $c$ since $\sum c'*c''\in
H_*(\sp{2}X,X)$.

To recapitulate, $\alpha+|c|$ is a cycle if, and only if, $\alpha$ is
a cycle and $c=0$. The only cycles in $C_*(W_2(X))$ are those that are
already cycles in the first summand $C_*(\sp{2}(X))$. On the other
hand, among these classes the only boundaries consist of boundaries in
$C_*(\sp{2}(X))$ and those of the form $\Delta_*(c)-c$ with $c$ a
cycle in $C_*(X)$ (in particular the only $0$-cycle is represented by
$x_0$). This proves our claim.
\end{proof}

\bre\label{pushoutsub} 
We could have noticed alternatively the existence
of a pushout diagram
$$\xymatrix{ X\vee
  X\ar[d]^{fold}\ar[r]^f&\sp{2}X\ar[d]^\alpha\\ X\ar[r]^{j_{x_0}\ \ \ \ }&
  \sub{3}(X,x_0) }$$ 
where $f(x,x_0) = [x,x]$ is the diagonal and $f(x_0,x) = [x,x_0]$.  We
can in fact deduce lemma \ref{computation} from this pushout. We can
also deduce that $\sub{3}(X,x_0)$ is simply connected if $X$ is. This
useful fact we use to establish proposition \ref{equivalent} next.
\ere

Note that lemma \ref{computation} above says that
$H_*(\sub{3}(X,x_0))$ only depends on $H_*(X)$ and on its coproduct
(i.e. on the cohomology of $X$). When $X$ is a suspension the
situation becomes simpler.  The following result is a nice
combination of lemmas \ref{null} and \ref{computation}.

\bpr\label{equivalent} There is a homotopy equivalence
$\sub{3}(\Sigma X,x_0)\simeq\bsp{2}(\Sigma X)$.\epr

\begin{proof} When $X$ is a suspension, all classes are primitive so that
$\Delta_*(c)=2c$ for all $c\in H_*(X)$. Combining Steenrod's
splitting (\ref{steenrod}) ;
$$H_*(\sp{2}X)\cong H_*(X)\oplus H_*(\sp{2}X,X)$$ with lemma
\ref{computation} we deduce immediately that $H_*(\sub{3}(\Sigma
X,x_0))\cong H_*(\bsp{2}(\Sigma X))$.  Both spaces are simply
connected (by remark \ref{pushoutsub} and theorem \ref{connectbsp})
and so it is enough to exhibit a map between them that induces this
homology isomorphism. Consider the map $\alpha : \sp{2}(\Sigma
X)\lrar\sub{3}(\Sigma X,x_0)$, $[x,y]\mapsto \{x,y,x_0\}$ as in the
proof of lemma \ref{computation}. Its restriction to $\Sigma X$ is
null-homotopic according to lemma \ref{null} and hence it factors
through the quotient $\bsp{2}(\Sigma X)\lrar\sub{3}(\Sigma
X,x_0)$. By inspection of the proof of lemma \ref{computation} 
we see that this map induces an isomorphism on homology.
\end{proof}

\bex\label{bsp2} A description of $\bsp{2}(S^k)$ is
given in (\cite{hatcher} , example 4K.5) from which we infer
that
$$\sub{3}(S^k,x_0)\simeq \Sigma^{k+1}\bbr P^{k-1}\ \ \  , \ \ k\geq 1$$ 
This generalizes the calculation in \cite{tuffley2} that
$\sub{3}(S^2,x_0)\simeq S^4$. \eex

\subsection{Homology Calculations}
We determine the homology of $\sub{3}(T,x_0)$ and $\sub{3}(T)$ where
$T$ is the torus $S^1\times S^1$. Symmetric products of surfaces are
studied in various places (see \cite{ks,tuffley2} and references
therein). Their homology is torsion free and hence particularly simple
to describe. We will write throughout $q: X^n\lrar\sp{n}X$ for the
quotient map and $q_*(a_1\tensor \ldots\tensor a_n) = a_1*a_2*\cdots
*a_n$ for its induced effect in homology (since our spaces are torsion
free we identify $H_*(X\times Y)$ with $H_*(X)\tensor H_*(Y)$).

\bco The inclusion $j: \sub{2}(T,x_0)\hookrightarrow\sub{3}(T,x_0)$ is
essential.\eco

\begin{proof} We will show that $j_*$ is non-trivial on
$H_2(\sub{2}(T,x_0)) = H_2(T)=\bbz$. Here $H_*(T)$ is generated by
$e_1,e_2$ in dimension one, and by the orientation class $[T]$ in
dimension two. The groups $H_*(\sp{2}T)$ are given as follows
\cite{ks} (the generators are indicated between brackets)
\begin{equation}\label{homtorus}
\tilde H_*(\sp{2}T) = \begin{cases} \bbz\{\gamma_2\},& \dim 4\\
\bbz\{e_1*[T], e_2*[T]\},& \dim 3\\
\bbz\{[T],e_1*e_2\},&\dim 2\\
\bbz\{e_1,e_2\},&\dim 1
\end{cases}
\end{equation}
where $\gamma_2$ is the orientation class $[\sp{2}T]$ ( $\sp{2}(T)$
is a compact complex surface).  Then $[T]*[T]=2\gamma_2$. Let
$\Delta$ be the diagonal into the symmetric square
$X\fract{\Delta}{\lrar} X\times X\fract{q}{\lrar}\sp{2}(X)$. Since
$\Delta_*([T]) = [T]\tensor 1 + e_1\tensor e_2 - e_2\tensor e_1 +
1\tensor [T]$, and since $q_*([T]\tensor 1) = q_*(1\tensor [T]) =
[T]$ and $q_*(e_1\tensor e_2) = -q_*(e_2\tensor e_1)= e_1*e_2$, we
see that
\begin{equation}\label{diagonal}
\Delta_*([T]) = 2[T] + 2e_1*e_2
\end{equation}
We can consider the composite
$$j_{x_0}: T\fract{\Delta}{\lrar}\sp{2}T\fract{\alpha}{\lrar}
\sub{3}(T,x_0)=\sp{2}T/\sim $$ where $\alpha$ is as in the proof of
lemma \ref{computation}. According to lemma \ref{computation}, using
the expression of the diagonal in (\ref{diagonal}), there are
classes $a = \alpha_*[T] , b = \alpha_*(e_1*e_2)$ with $a = -2b\neq
0$. But $(j_{x_0})_*[T]=(\alpha\circ \Delta)_*[T] =
\alpha_*([T])=a,$ and this is non-zero as desired.
\end{proof}

\bre We can of course complete the calculation of
$H_*(\sub{3}(T,x_0))$ from lemma \ref{computation}. Under
$\alpha_*$, $e_i\mapsto 0$ (primitive classes map to $0$),
$e_1*e_2\mapsto b$, $[T]\mapsto a = -2b$, $e_i*[T]\mapsto c_i,$ and
$\gamma_2\mapsto d,$ so that
$$ H_1 = 0 \ ,\ H_2 = \bbz\{a\}\ ,\ H_3=\bbz\{c_1,c_2\} \ ,\ H_4 =
\bbz\{d\}
$$
It is equally easy to write down the homology groups for
$\sub{3}(S,x_0)$ for any genus $g\geq 1$ surface, orientable or not.
\ere

Next we analyze the inclusion $T\hookrightarrow\sub{3}T$ in the case
of the torus (compare \cite{tuffley2}). The starting point is the
pushout (\ref{pushout3}) and the associated Mayer-Vietoris sequence
$$\cdots\lrar H_*(T\times T)\fract{q_*\oplus i_*}{\ra 3}
H_*(\sp{2}T)\oplus H_*(\sp{3}T)\fract{g_* - \pi_*}{\ra 3}
H_*(\sub{3}T)\lrar H_{*-1}(T\times T)\lrar\cdots
$$
where $q: T\times T\lrar\sp{2}T$ is the quotient map, $i(x,y)=x^2y$,
$g : \sp{2}T\hookrightarrow\sub{3}T$ is the inclusion (here we have
identified $\sp{2}T$ with $\sub{2}T$) and $\pi :
\sp{3}T\lrar\sub{3}T$ is the projection. We focus on degree $2$ and
follow \cite{ks} for the next computations.

We have $H_2(T\times T) = \bbz^2$ generated by $[T]\tensor 1$ and
$1\tensor [T]$, $H_2(\sp{2}T) = \bbz^2=H_2(\sp{3}T)$ generated by a
class of the same name $[T] = q_*([T]\tensor 1)=q_*(1\tensor [T])$
and by $e_1*e_2$ ; see (\ref{homtorus}). To describe the effect of
$i_*$ we write it as a composite
$$i : T\times T\fract{\Delta\times 1}{\ra 3} T\times T\times T
\fract{q}{\lrar} \sp{3}T$$ This gives  $i_*([T]\tensor 1) = 2[T] +
2e_1*e_2$ as in (\ref{diagonal}), while $i_*(1\tensor [T])=[T]$. The
Mayer-Vietoris then looks like
\begin{eqnarray*}
\cdots\lrar\ \bbz^2 &\fract{q_*\oplus i_*}{\ra 3}&
\bbz^2\oplus\bbz^2\fract{g_* - \pi_*}{\ra 3} H_2(\sub{3}T)\lrar
H_{1}(T\times T)\lrar\cdots\\
(1,0)&{\longmapsto} &((1,0), (2,2))\\
(0,1)&\longmapsto &((1,0), (1,0))
\end{eqnarray*}
This sequence is exact. Observe that the class $((2,2),(0,0))$ is
not in the kernel of $g_*-\pi_*$  because it cannot be in the image
of $q_*\oplus i_*$. This means that $g_*(2,2) \neq 0$. This is all
we need to derive the non-nullity of the map $j:
X\hookrightarrow\sub{3}X$.

\bco $j_*([T])\neq 0$. \eco

\begin{proof} The inclusion $j$ is the composite
$$j : X\fract{\Delta}{\lrar}X\times X\fract{\pi}{\lrar}
\sp{2}X\fract{g}{\lrar}\sub{3}X$$ so that $j_*([T]) = g_*(2,2),$ and
this is non-trivial as asserted above.
\end{proof}


\section{The Top Dimension}\label{topdim}

Using facts about orientability of configuration spaces of closed
manifolds (\cite{braid} for example) we slightly elaborate on
\cite{handel} and (\cite{tuffley2} theorem 3).

\bpr\label{top} Suppose $M$ is a closed manifold of dimension $d\geq
2$. Then
$$
H_{nd}(\sp{n}M;\bbz ) =
\begin{cases} \bbz &\hbox{if}\ d\ \hbox{even and $M$ orientable}\\
0&\hbox{if}\ d\ \hbox{odd or $M$ non-orientable}
\end{cases}
$$
For mod-$2$ coefficients, $H_{nd}(\sp{n}M;\bbf_2 ) = \bbf_2$. In all
cases the map
$$H_{nd}(\sp{n}M )\lrar H_{nd}(\sub{n}M )$$
is an isomorphism (Corollary \ref{same}). \epr

\begin{proof}
When $d=2$ the claim is immediate since, as is well known, $\sp{n}M$
is a closed manifold (orientable if and only if $M$ is; see
\cite{wagner}). Generally our statement follows from the fact that
$\sp{n}(X)$ is an orbifold with codimension $>1$ singularities, and
hence its top homology group is that of a manifold. More explicitly
in our case, let's denote by $B(M,n)$ the configuration space of
finite sets of cardinality $n$ in $M$; that is
$$B(M,n) = \sp{n}M - {\bf\Delta}^n = \sub{n}M - \sub{n-1}M$$
where ${\bf\Delta}^n$ is the singular set consisting of tuples with
at least one repeated entry (the image of the fat diagonal as
defined in \S\ref{basic}).  By Poincar\'e duality suitably applied
(\cite{braid}, lemma 3.5)
\begin{equation}\label{dual}
H^i(B(M,n);\pm \bbz)\cong H_{nd-i}(\sp{n}M, {\bf\Delta}^n;\bbz )
\end{equation}
where $\pm\bbz$ is the orientation sheaf.  By definition
$$H^i(B(M,n),\pm\bbz ) = H^i(Hom_{Br_n(M)}(C_*(\tilde B(M,n)),\bbz ))$$
where $Br_n(M)=\pi_1(B(M,n))$ is the braid group of $M$, $\tilde
B(M,n)$ is the universal cover of $B(M,n)$ and the action of the
class of a loop on $\bbz$ is multiplication by $\pm 1$ according to
whether the loop preserves or reverses orientation.  It is known
that $B(M,n)$ is orientable if and only if $M$ is orientable and
even dimensional (\cite{braid}, lemma 2.6). That is we can replace
$\pm\bbz$ by $\bbz$ if $M$ is orientable and $d$ is even.

Since ${\bf\Delta}^n$ is a subcomplex of codimension $d$ in
$\sp{n}M$, we have $H_{nd-i}(\sp{n}M, {\bf\Delta}^n )\cong
H_{nd-i}(\sp{n}M)$ for $i < d-1$ . In particular, for $i=0$ we
obtain
\begin{equation}\label{dualcn}
H^0(B(M,n);\pm\bbz )\cong H_{nd}(\sp{n}M;\bbz )
\end{equation}
If $M$ is even dimensional and orientable, $H^0(B(M,n); \pm\bbz )
\cong H^0(B(M,n); \bbz )=\bbz$ since $B(M,n)$ is connected if $\dim
M\geq 2$. If $\dim M$ is odd or $M$ is non-orientable, then $B(M,n)$
is not orientable and $H^0(B(M,n);\pm\bbz )=0$ (this is because
$H^0(B(M,n);\pm\bbz )$ is the subgroup $\{m\in\bbz\ |\ gm = m,\
\forall g\in\bbz[\pi_1(B(M,n)]\}$). This establishes the claim for
the symmetric products and hence for the finite subset spaces
according to corollary \ref{same}.
\end{proof}

\bex For $k\geq 2$ we have
$H_{2k}(\sp{2}S^k)=H_{2k}(\bsp{2}S^k)=H_{k-1}(\bbr P^{k-1})$ (see
example \ref{bsp2}) and this is $\bbz$ or $0$ depending on whether
$k$ is even or odd as predicted by proposition \ref{top}. \eex

\subsection{The Case of the Circle} When $M=S^1$ , proposition
\ref{top} is not true anymore since $\sp{n}S^1\simeq S^1$ for all
$n\geq 1$, while $\sub{n}(S^1)$ is either $S^n$ or $S^{n-1}$
depending on whether $n$ is odd or even \cite{mp, tuffley1}. It is
still possible to describe in this case the quotient map
$\sp{n}(S^1)\lrar\sub{n}(S^1)$ explicitly.

A beautiful theorem of Morton asserts that the multiplication map
$$\sp{n+1}(S^1)\lrar S^1$$
is an $n$-disc bundle $\eta_n$ over $S^1$ which is orientable if,
and only if, $n$ is even \cite{morton}. A close scrutiny of Morton's
proof shows that the sphere bundle associated to $\eta_n$ consists
of the image of the fat diagonal ${\bf\Delta}^{n+1}$; i.e. the
singular set. If $Th(\eta_n)$ is the Thom space of $\eta_n$, then
\begin{equation}\label{one}
Th(\eta_n) = \sp{n+1}(S^1)/{\bf\Delta}^{n+1} =
\sub{n+1}S^1/\sub{n}S^1
\end{equation}
Since $\eta_n$ is trivial when $n=2k$ is even, it follows that
\begin{equation}\label{two}
Th(\eta_{2k}) = S^{2k}\wedge S^1_+ = S^{2k+1}\vee S^{2k}
\end{equation}

But as pointed out above, $ \sub{2k+1}(S^1)\simeq S^{2k+1}$. The map
$\sp{2k+1}(S^1)\lrar\sub{2k+1}(S^1)$ factors  through the Thom space
(\ref{two}) and the top cell maps to the top cell. Combining
(\ref{one}) and (\ref{two}) it is immediate to see that

\ble The map $Th(\eta_{2k})\lrar\sub{2k+1}(S^1)$ restricted to the
first wedge summand in (\ref{two}) induces a map
$S^{2k+1}\lrar\sub{2k+1}(S^1)$ which is a homotopy equivalence.\ele


\section{Manifold Structure}\label{manifold}

In this last section we prove Theorem \ref{main3}. We distinguish
three cases : when the dimension of the manifold is $d>2$, $d=2$ or
$d=1$.

\ble\label{n3d2} Suppose $X$ is a manifold of dimension $d>2.$  Then
$Sub_nX$ is never a manifold if $n\ge 2.$\ele

\begin{proof}
Consider the projection $X^n\lrar\sub{n}X$ given by identifying
tuples whose sets of coordinates are the same. This projection
restricts to an $n!$ regular covering between the complements $\pi_n
: X^n-\Lambda^n\lrar \sub{n}X - \sub{n-1}X,$ where $\Lambda^n$
as before is the fat diagonal in $X^n$. Suppose $\sub{n}X$ is a manifold of
dimension $nd$ (necessarily). Pick a point in $\sub{n-1}X$ and an
open chart $U$ around it. Now $U\cong\bbr^{nd}$ and
$Y=U\cap\sub{n-1}X$ is a closed subset in $U$. We can apply
Alexander duality to the pair $(Y,U)$ and obtain
$$H_{nd-i-1}(U-Y)\cong H^i(Y)$$
But $Y\subset\sub{n-1}(X)$ is an open subspace in a simplical
complex of dimension $(n-1)d$; therefore $H^{nd-2}(Y)=0$ (since $d>
2$) and so $H_1(U-Y)=0$. We can now use an elementary observation of
Mostovoy \cite{mostovoy} to the effect that since $U-Y$ is covered
by $\pi_n^{-1}(U-Y),$ a connected \'etale cover of degree $n!$, then
it is impossible for $H_1(U-Y)$ to be trivial since the monodromy
gives a surjection $\pi_1(U-Y)\lrar\sn,$ and hence
  a non-trivial map $H_1(U-Y)\lrar\bbz_2$.
\end{proof}

Theorem 2.4 of \cite{wagner} shows that our Lemma \ref{n3d2} is
valid if $d=2$ and $n>2$ as well. As opposed to the geometric
approach of Wagner, we provide below a short homological proof of
this result.

\ble\label{d2n2} Suppose $X$ is a closed topological surface. Then
$\sub{n}X$ is a manifold if and only if $n=2$. \ele

\begin{proof}
We will show that if $n\geq 3$, then $\sub{n}(X)$ cannot even have
the homotopy type of a closed manifold by showing that it doesn't
satisfy Poincar\'e duality. We rely on results of \cite{ks} that
give a simple description of a CW decomposition of a space
$\wsp{n}X$ homotopy equivalent to $\sp{n}X$ when $X$ is a two
dimensional complex. Since $X$ is a closed two dimensional manifold,
it has a cell structure of the form $X = \bigvee^r S^1\cup D^2$
where $D^2$ is a two dimensional cell attached to a bouquet of
circles. Each circle corresponds in the cellular chain complex for
$\wsp{n}X$ to a one-dimensional cell generator $e_i$, $1\leq i\leq
r$, while the two dimensional cell is represented by $D$. This chain
complex has a concatenation product $*: C_*(\wsp{r}X)\tensor
C_*(\wsp{s}X)\lrar C_*(\wsp{r+s}X)$ under which these cells map to
product cells. The full cell complex for $\wsp{n}X$ is made up of
all products of the form
$$e_{i_1}*\cdots * e_{i_\ell}*\sp{k}D\ \ \ ,\ \ \ i_1+\cdots +i_\ell +
k\leq n$$ where $i_r\neq i_s$ if $r\neq s$, and where $\sp{k}D$ is a
$2k$-dimensional cell represented geometrically by the $k$-th
symmetric product of $D^2$. The boundary $\partial$ is a derivation
and is completely determined on generators by $\partial e_i = 0$ and
$\partial \sp{n}D =
\partial D*\sp{n-1}D$.

If $X = \bigvee^r S^1\cup D$ is a closed manifold, then in mod-$2$
homology, $\partial D = 0$ (the top cell). This implies of course
that $\partial\sp{n}D = 0$ (the top cell of $\sp{n}X$), while
$H_{2n-1}(\sp{n}X,\bbz_2)\cong\bbz_2^r$ with generators $e_i*\sp{n-1}D$.
This shows in particular that $H_{2n-1}(\sp{n}X;\bbz_2)\neq 0$ if
$r\geq 1$; that is if $X$ is not the two sphere. Observe that this
calculation is compatible with Theorem 2 of \cite{tuffley2}.

Now we know that $\sub{n}X$ is simply connected if $n\geq 3$.
Suppose $\sub{n}X$ is a closed manifold, then by Poincar\'e duality
$H_{2n-1}(\sub{n}X;\bbz_2) = H_1(\sub{n}X;\bbz_2) = 0$. But recall
the pushout diagram (\ref{pushout2}) and its associated
Mayer-Vietoris exact sequence
$$H_{2n-1}(\Delta_n)\lrar H_{2n-1}(\sub{n-1}X)\oplus H_{2n-1}(\sp{n}X)
\lrar H_{2n-1}(\sub{n}X)\lrar H_{2n-2}(\Delta_n)\lrar\cdots
$$
Since $\Delta_n$ and $\sub{n-1}X$ are $(2n-2)$-dimensional
subcomplexes of $\sub{n}X$, their homology in degree $2n-1$
vanishes. The sequence above becomes
$$0\lrar H_{2n-1}(\sp{n}X)\lrar H_{2n-1}(\sub{n}X)\lrar
H_{2n-2}(\Delta_n)\lrar\cdots$$ and $H_{2n-1}(\sp{n}X)$ injects into
$H_{2n-1}(\sub{n}X)$. When $H_1(X)\neq 0$; that is when $X$ is not
the sphere, $H_{2n-1}(\sub{n}X)$ is non-trivial
contradicting Poincar\'e duality.

We are left with the case $\sub{n}(S^2)$ and $n\geq 3$.  Here we have
to rely on a calculation of Tuffley \cite{tuffley2} who shows that
\begin{equation}\label{calcul} H_{2n-2}(\sub{n}(S^2)) =
\bbz\oplus\bbz_{n-1} \end{equation} But $\sub{n}(S^2)$ is
$2$-connected according to Theorem \ref{main1} and Poincar\'e
duality is violated in this case as well.
\end{proof}

\bre A computation of the homology of $\sub{n}(S^2)$ for all $n$ and
various field coefficients will appear in \cite{sadok}. It is
however straightforward using the Mayer-Vietoris sequence for the
pushout (\ref{pushout3}) to show that
\begin{equation} \tilde H_*(\sub{3}S^2)\cong\begin{cases}\bbz&, *=6\\
  \bbz\oplus\bbz_2&, *=4
\end{cases}
  \end{equation}
Similar computations appear in \cite{biro, tuffley2, taamallah}. \ere

Finally we address the case $d=1$. Write $I=[0,1], \dot{I}=(0,1)$.
First of all $\sp{n}(I)\cong I^n$. In fact this is precisely the
$n$-simplex since any point of $\sp{n}(I)$ can be written uniquely
as an $n$-tuple $(x_1,\ldots, x_n)$ with $0\leq x_1\leq\cdots\leq
x_n\leq 1$. The quotient map $q_2 : \sp{2}(I)\lrar\sub{2}(I)$ is a
homeomorphism and hence every interior point of $\sub{2}(I)$ has a
manifold neighborhood.  The same for $n=3$ since $\sp{3}(I)$ is the
three simplex
$$\{(x_1,x_2,x_3)\ |\ 0\leq x_1\leq x_2\leq x_3\leq 1\}$$ with $4$
faces: $F_1: \{x_1=0\}$, $F_2: \{x_1=x_2\}$, $F_3:\{x_2=x_3\}$ and
$F_4:\{x_3=1\}$, and the quotient map $q_3 :\sp{3}(I)\rightarrow\sub{3}(I)$
identifies the faces $F_2$ and $F_3$. Such an identification gives
again $I^3$ and $\sub{3}(\dot{I})$ is this simplex
with two faces removed \cite{rose}.  For $n>3$, the corresponding map
$q_n$ identifies various faces of the simplex $\sp{n}(I)$ to obtain
$\sub{n}(I),$ but this fails to give a manifold structure on the
quotient for there are just too many ``branches" that come together at
a single point in the image of the boundary of this simplex. This is
made precise below.

\ble\label{subs1} $\sub{n}(S^1)$ is a closed manifold if and only if
$n=1,3$. \ele

Observe that if $n$ is even, then $\sub{n}S^1$ cannot be a closed
manifold for a simple reason: no closed manifold of dimension $n$
can be homotopic to a sphere of dimension $n-1$.

\begin{proof} (of Lemma \ref{subs1} following \cite{wagner},  Theorem
2.3). Let $M$ be a manifold and $D$ a disc neighborhood of a point
$x\in M$. Then an open neighborhood of $x\in\sub{n}(M)$ is
$\sub{n}(D)$. So if $\sub{n}(D)$ is not a manifold, then neither is
$\sub{n}(M)$. To prove lemma \ref{subs1} we will argue as in
\cite{wagner} that $\sub{n}(\bbr)$ is not a manifold for $n\geq 4$.

For a metric space $X$ (with metric $d$),  non-empty subsets
$S,T\subset X,$ and fixed elements $s\in S,t\in T,$ we define
\begin{eqnarray*}
d(s,T)&=&\inf\{d(s,t) \bigm | t\in T\}\\
d(S,t)&=&\inf\{d(s,t) \bigm | s\in S\}
\end{eqnarray*}
Then the Hausdorff metric $D$ on $\sub{n}(X)$ is defined to be
$$D(S,T):= \hbox{sup}\{d(s,T) ,d(t,S)\ | \ s\in S, t\in T\}$$
Thus $D(S,T)<\epsilon$ means that each $s\in S$ is within an
$\epsilon$-neighborhood of some point in $T$ and each $t\in T$ is
within an $\epsilon$-neighborhood of some point in $S$.

We wish to show that $\sub{n}(\bbr )$ for $n\geq 4$ is not
homemorphic to $\bbr^n$. Pick $S=\{1,2,\ldots, n-1\}$ in
$\sub{n-1}(\bbr )$ and for each $i$ consider the open set $C_i$ (in
the Hausdorff metric) of all subsets $\{p_1,\ldots,
p_{n-1},q_i\}\in\sub{n}(\bbr )$ such that $p_j\in (j-\frac{1}{2},
j+\frac{1}{2})$ and $q_i\in (i-\frac{1}{2}, i+\frac{1}{2})$. We then
see that $C_i$ is the subset with one or two points in the
$\frac{1}{2}$-neighborhood of $i$ and a single point in the
$\frac{1}{ 2}$-neighborhood of $j$ for $i\neq j$. Note that
$C_i\subset U$ where $U=\{T\in\sub{n}(\bbr )\ |\ D(S,T)<1/2\}$.
Observe that
\begin{eqnarray*}
C_1 &=& \sub{2}\left(\frac{1}{2}, \frac{3}{2}\right)\times
\left(\frac{3}{2}, \frac{5}{2}\right)\times\cdots\times
\left(n-1-\frac{1}{2}, n-1 + \frac{1}{2}\right)
\end{eqnarray*}
This is an $n$-dimensional manifold with boundary
$V=U\cap\sub{n-1}(\bbr )$ and in fact one has
$$C_i = \left\{T\in U : T\cap \left(i-\frac{1}{2},i+\frac{1}{2}\right)\
\hbox{has $1$ or $2$ points}\right\} \cup V$$ Clearly $C_1\cup
C_2\cup\cdots \cup C_{n-1} = U$ and more importantly all these open
sets have a common boundary at $V$; i.e. $C_i\cap C_j = V$.  If
$n\geq 4$, we can choose at least three such $C_i$; say
$C_1,C_2,C_3$. Then $C_1\cup C_2$ is an open $n$-dimensional
manifold (union over the common boundary $V$). It must be contained
in the interior of $\sub{n}(\bbr )$ and hence must be open there if
$\sub{n}(\bbr )$ were to be an $n$-dimensional manifold. But
$C_1\cup C_2$ is not open in $\sub{n}(\bbr )$ since every
neighborhood of $\{1,2,\ldots, n-1\}$ must meet $C_3-V$ which is
disjoint from $C_1\cup C_2$ (i.e. ``too many'' branches come
together at that point).
\end{proof}

We conclude this paper with the following cute theorem of Bott,
which is the most significant early result on the subject.

\bco (Bott) There is a homeomorphism $\sub{3}(S^1)\cong S^3$. \eco

\begin{proof} It has been known since Seifert that the
Poincar\'e conjecture holds for Seifert manifolds; that is if a
Seifert $3$-manifold is simply connected  then it is homeomorphic to
$S^3$\footnote{We thank Peter Zvengrowski for reminding us this fact}.
Clearly $\sub{3}(S^1)$ is a Seifert manifold where the action of
$S^1$ on a subset is by multiplication on elements of that subset.
Since it is simply connected (corollary \ref{fundamental}), the
claim follows. Note that the $S^1$-action has two exceptional fibers
consisting of the orbits of $\{1,-1\}$ and $\{1,j,j^2\}$ where
$j=e^{2\pi i/3}$ (compare \cite{tuffley1}).
\end{proof}


\addcontentsline{toc}{section}{Bibliography}
\bibliography{biblio}
\bibliographystyle{plain[10pt]}

\end{document}